\documentclass[12pt]{article}
\usepackage{amssymb}
\begin{document}
\renewcommand{\thefootnote}{\fnsymbol{footnote}}
\newpage
\pagestyle{empty}
\setcounter{page}{0}
\newcommand{\dotimes}{\stackrel{}{\dot{\otimes}}}
\renewcommand{\thesection}{\arabic{section}}
\renewcommand{\theequation}{\thesection.\arabic{equation}}
\newcommand{\sect}[1]{\setcounter{equation}{0}\section{#1}}
\newfont{\twelvemsb}{msbm10 scaled\magstep1}
\newfont{\eightmsb}{msbm8}
\newfont{\sixmsb}{msbm6}
\newfam\msbfam
\textfont\msbfam=\twelvemsb
\scriptfont\msbfam=\eightmsb
\scriptscriptfont\msbfam=\sixmsb
\catcode`\@=11
\def\Bbb{\ifmmode\let\next\Bbb@\else
\def\next{\errmessage{Use \string\Bbb\space only in math mode}}\fi\next}
\def\Bbb@#1{{\Bbb@@{#1}}}
\def\Bbb@@#1{\fam\msbfam#1}
\newfont{\twelvegoth}{eufm10 scaled\magstep1}
\newfont{\tengoth}{eufm10}
\newfont{\eightgoth}{eufm8}
\newfont{\sixgoth}{eufm6}
\newfam\gothfam
\textfont\gothfam=\twelvegoth
\scriptfont\gothfam=\eightgoth
\scriptscriptfont\gothfam=\sixgoth
\def\frak{\frak@}
\def\frak@#1{{\fam\gothfam{{#1}}}}
\def\frak@@#1{\fam\gothfam#1}
\catcode`@=12
%
%
%
\def\CC{{\Bbb C}}
\def\NN{{\Bbb N}}
\def\QQ{{\Bbb Q}}
\def\RR{{\Bbb R}}
\def\ZZ{{\Bbb Z}}
\def\cA{{\cal A}}          \def\cB{{\cal B}}          \def\cC{{\cal C}}
\def\cD{{\cal D}}          \def\cE{{\cal E}}          \def\cF{{\cal F}}
\def\cG{{\cal G}}          \def\cH{{\cal H}}          \def\cI{{\cal I}}
\def\cJ{{\cal J}}          \def\cK{{\cal K}}          \def\cL{{\cal L}}
\def\cM{{\cal M}}          \def\cN{{\cal N}}          \def\cO{{\cal O}}
\def\cP{{\cal P}}          \def\cQ{{\cal Q}}          \def\cR{{\cal R}}
\def\cS{{\cal S}}          \def\cT{{\cal T}}          \def\cU{{\cal U}}
\def\cV{{\cal V}}          \def\cW{{\cal W}}          \def\cX{{\cal X}}
\def\cY{{\cal Y}}          \def\cZ{{\cal Z}}
\def\qed{\hfill \rule{5pt}{5pt}}
\def\id{\mbox{id}}
\def\ggo{{\frak g}_{\bar 0}}
\def\uqggo{\cU_q({\frak g}_{\bar 0})}
\def\uqggp{\cU_q({\frak g}_+)}
\def\half{\frac{1}{2}}
\def\btf{\bigtriangleup}
\newtheorem{lemma}{Lemma}
\newtheorem{prop}{Proposition}
\newtheorem{theo}{Theorem}
\newtheorem{Defi}{Definition}

\vfill
\vfill
\begin{center}

{\LARGE {\bf {\sf
Entanglement via Barut-Girardello coherent state for $su_{q}(1,1)$
quantum algebra: bipartite composite system
}}} \\[0.8cm]

{\large
R. Chakrabarti\footnote{E-mail:
{\tt ranabir@imsc.res.in}} and S. S. Vasan
}

\begin{center}
{\em
Department of Theoretical Physics, University of Madras,
Guindy Campus,\\
Chennai, 600 025, India. \\[2.8cm]
}
\end{center}

\end{center}

\smallskip

\smallskip

\smallskip

\smallskip

\smallskip

\smallskip

\begin{abstract}
\noindent Using noncocommutative coproduct properties of the
quantum algebras, we introduce and obtain, in a bipartite
composite system, the Barut-Girardello coherent state for the
$q$-deformed $su_{q}(1,1)$ algebra. The quantum coproduct
structure ensures this normalizable coherent state to be entangled.
The entanglement disappears in the classical $q \rightarrow 1$
limit, giving rise to a factorizable state.
\end{abstract}

\vfill
\newpage

\pagestyle{plain}

\sect{Introduction}
Entanglement is the key distinguishing feature of quantum
mechanics setting it apart from classical physics. A quantum state
of a composite system, consisting of two or more subsystems, is
entangled if it cannot be factorized into direct product of the
states of the subsystems. Entangled states are useful in quantum
information processing such as  quantum teleportation
\cite{BBCJPW93}, quantum key distribution \cite{E91} and superdense
coding \cite{BW92}. Studying quantum information theory using
entangled coherent states has recently received much attention
\cite{S92}-\cite{MMS00}. In a related context the coherent states
of the $su(2)$ and the $su(1, 1)$ algebras were studied
\cite{WSP00}.

\par

The purpose of the present work is to extend the horizon of studies
on entangled nonorthogonal states so as to incorporate systems with
quantum algebraic symmetries \cite{KS97}. Composite systems with
quantum symmetries, such as anyons for instance \cite{CGM93},
are natural candidates for studying entangled states. The
reason for this lies in the noncocommutativity of the coproduct map
of the generators of the quantum algebras. As a demonstration of
this property we here analytically obtain, in a bipartite composite
system, the Barut-Girardello coherent state \cite{BG71} for the
$su_{q}(1, 1)$ quantum algebra \cite{KS97}. The entangled coherent
state (\ref{eq:qbBGcost}) obtained here is not factorizable in the
quantum states of its subsystems for a generic value of the
deformation parameter $q$. As $q \rightarrow 1$ in the classical
limit, the entanglement in the state (\ref{eq:qbBGcost}) disappears
reducing it to the factorized classical form (\ref{eq:clfact}). For
the purpose of setting the framework we first study, in the context
of single-node systems, the Barut-Girardello coherent state for a
general class of deformed $su(1, 1)$ algebras. In particular, we
explicitly demonstrate completeness relation for the $q$-deformed
$su_{q}(1, 1)$ Barut-Girardello coherent states in terms of an
ordinary integral over the complex plane. Entangled coherent state
in a bipartite composite system is studied in Sec. 3.

\sect{Barut-Girardello coherent states for the deformed $su (1,1)$
algebras: single-node systems}
The generators $( K_{0}, K_{\pm} )$ of the classical $su (1,1)$
algebra satisfy the defining commutation relations
\begin{equation}
[ K_{0}, K_{\pm} ] = \pm K_{\pm},\qquad
[ K_{-}, K_{+} ] = 2 K_{0}
\label{eq:clalg}
\end{equation}
and maintain the hermiticity constraints
$( K_{0}^{\dagger} = K_{0}, K_{+}^{\dagger} = K_{-} )$. The Casimir
element of the algebra reads
\begin{equation}
C = K_{0}^{2} - K_{0} - K_{+} K_{-}.
\label{eq:clcas}
\end{equation}
For the discrete series of representations the basis states read
$\{| n,k \rangle\,| n = 0, 1, 2,\cdots ; 2 k = \pm 1, \pm 2,
\cdots \}$ and the irreducible representations are parametrized by
a single number $k$: $C = k ( k - 1 ) {\sf I}$. An arbitrary
irreducible representation reads
\begin{eqnarray}
K_{0}\,| n, k \rangle &=& (n + k)\,| n, k \rangle,\nonumber\\
K_{+}\,| n, k \rangle &=& \sqrt {(n + 1)\,(n + 2 k)}
\,| n + 1, k \rangle,\nonumber\\
K_{-}\,| n, k \rangle &=& \sqrt {n\,(n + 2 k - 1)}
\,| n - 1, k \rangle.
\label{eq:clrep}
\end{eqnarray}
We assume that the set of states described above form a complete
orthonormal basis. The primitive coproduct structure of the
classical generators is given by
\begin{equation}
\btf(K_{i}) = K_{i} \otimes 1 + 1 \otimes K_{i}\qquad\forall
i \in (0, \pm).
\label{eq:clcoprod}
\end{equation}

\par

General nonlinear deformations of the $su (2)$ and the $su (1,1)$
algebras were considered in  \cite{P90}, \cite{R91} and
\cite{DQ93}. In particular it was observed in \cite{DQ93} that for
a class of these nonlinear deformations exponential spectra occur
in the carrier space of the unitary irreducible representations.
In a parallel development in the context of oscillators Manko
{\it et al.} \cite{MMSZ97} introduced, {\it via} nonlinear maps,
the notion of $f$-oscillators as a generalization of the standard
$q$-oscillators. Using the generalized deformed oscillator algebra
they also constructed the nonlinear $f$-coherent states. We
follow their approach and construct nonlinear Barut-Girardello
coherent states for a generalized deformed $su (1,1)$ algebra. We
review the technique here {\it in extenso} as our future
construction of a bipartite Barut-Girardello coherent state for the
$q$-deformed $su_{q}(1, 1)$ algebra involves similar methodology.
The generators of the nonlinear algebra are introduced {\it via} an
invertible map on the corresponding classical generators:
\begin{equation}
{\sf K}_{0} = K_{0},\quad
{\sf K}_{+} = f (K_{0})\,K_{+},\quad
{\sf K}_{-} = K_{-}\,f (K_{0}),
\label{eq:fsumap}
\end{equation}
where $f (K_{0})$ is an arbitrary operator-valued real function.
The complete orthonormal states introduced in (\ref{eq:clrep})
also constitute a carrier space of the deformed generators. The
form of the mapping function $f (K_{0})$ determines whether the
realization is irreducible or not. The generalized deformed
generators follow a nonlinear algebra:
\begin{equation}
[ {\sf K}_{0}, {\sf K}_{\pm} ] = \pm {\sf K}_{\pm},\qquad
[ {\sf K}_{-}, {\sf K}_{+} ] = {\sf F} ({\sf K}_{0}),
\label{eq:fsualg}
\end{equation}
where ${\sf F} ({\sf K}_{0}) =
({\sf K}_{0} + k)\,({\sf K}_{0} - k + 1)\,(f({\sf K}_{0} + 1))^{2}
- ({\sf K}_{0} - k)\,({\sf K}_{0} + k - 1)\,(f({\sf K}_{0}))^{2}$.
As we are concerned with a single-node system in this section,
here we do not consider the coalgebraic properties of the above
$f$-deformed $su_{f} (1,1)$ algebra. A suitable induced coproduct
structure may be realized for the $su_{f} (1,1)$ algebra. From the
point of view of the classical algebra, the deformed generators
introduced in (\ref{eq:fsumap}) may be regarded as nonlinear
operators which may be of significance in a particular physical
situation.

\par

For a single-node system the nonlinear Barut-Girardello coherent
state for the deformed algebra (\ref{eq:fsualg}) is defined as an
eigenstate of the generator ${\sf K}_{-}$:
\begin{equation}
{\sf K}_{-}\,|\alpha, k\rangle_{f} = \alpha\,|\alpha, k\rangle_{f},
\qquad \alpha \in {\mathbb C}.
\label{eq:bgcsn}
\end{equation}
The above coherent state may be expanded in terms of the basis
states introduced in (\ref{eq:clrep}):
\begin{equation}
|\alpha, k\rangle_{f} = \sum_{n = 0}^{\infty} c_{n}^{(f)}\,
|n, k\rangle.
\label{eq:snexp}
\end{equation}
Inserting the expansion (\ref{eq:snexp}) in the defining
relation (\ref{eq:bgcsn}), we, {\it via} the use of the map
(\ref{eq:fsumap}), obtain the recurrence relation
\begin{equation}
c_{n + 1}^{(f)} = \frac{\alpha}{f (n + k +1)\,\sqrt{(n + 1)\,
(n + 2 k)}}\,c_{n}^{(f)},
\label{eq:snrec}
\end{equation}
whose solution reads
\begin{equation}
c_{n}^{(f)} = N_{f}\,\frac{\alpha^{n}}{[f (n + k)]!\,\sqrt{n!\,
\Gamma (n + 2 k)}},\qquad
[f (n + k)]! = \prod_{j = 1}^{n}\,f (j + k).
\label{eq:snrecsol}
\end{equation}
The normalization condition $_{f}\langle\alpha, k|
\alpha, k\rangle_{f} = 1$  fixes the constant $N_{f}$:
\begin{equation}
N_{f}^{-2} = \sum_{n = 0}^{\infty}\,\frac{|\alpha|^{2 n}}
{\left([f (n + k)]!\right)^{2}\,n!\,\Gamma (n + 2 k)}.
\label{eq:Nf}
\end{equation}
The preceding derivation yields the normalized nonlinear
Barut-Girardello coherent state for the $f$-deformed $su_{f}(1,1)$
algebra:
\begin{equation}
|\alpha, k\rangle = N_{f}\,\sum_{n = 0}^{\infty} \frac{\alpha^{n}}
{[f (n + k)]!\,\sqrt{n!\,\Gamma (n + 2 k)}}\,|n, k\rangle.
\label{eq:snbgcs}
\end{equation}
Using the properties of the carrier space (\ref{eq:clrep}) the
single-node $f$-coherent state obtained above may be expressed in
terms operator-valued hypergeometric function as follows:
\begin{eqnarray}
|\alpha, k\rangle_{f} &=& c_{0}^{(f)}\exp \left(\alpha\,
(f({\sf K}_{0}))^{-2}\,{\sf K}_{+}\,({\sf K}_{0} + k)^{-1}\right)
\,|0, k\rangle\nonumber\\
\phantom{|\alpha, k\angle_{f}} &=& c_{0}^{(f)}\,
_{0}F_{1}\left(\underline{\phantom{a}};2 k;\alpha\,
(f({\sf K}_{0}))^{-2}\,{\sf K}_{+}\right)\,|0, k\rangle.
\label{eq:cshyp}
\end{eqnarray}
In the classical limit $f(K_{0}) \rightarrow 1$, the $f$-coherent
state constructed in (\ref{eq:snbgcs}) reduces to the
Barut-Girardello coherent state $|\alpha, k\rangle$ for the
classical $su(1,1)$ algebra \cite{BG71}:
\begin{equation}
|\alpha, k\rangle_{f}\,\longrightarrow\,|\alpha, k\rangle =
\frac{|\alpha|^{k - \half}} {\sqrt{I_{2 k - 1}(2 |\alpha|)}}\,
\sum_{n = 0}^{\infty}\,\frac{\alpha^{n}}{\sqrt{n!\,
\Gamma (n + 2 k)}}\,|n, k\rangle,
\label{eq:cslim}
\end{equation}
where the modified Bessel function of the first kind is given by
\begin{equation}
I_{m}(2 z) = \sum_{n = 0}^{\infty}\frac{z^{m + 2 n}}{n!\,
\Gamma (m + n + 1)}.
\label{eq:modbf}
\end{equation}
Parallel to its classical analog, the set of nonlinear coherent
states $|\alpha, k\rangle_{f}$ exhibits the important property of
completeness (actually overcompleteness). Using the polar
decomposition $\alpha = \rho\,\exp (i\,\theta)$ and integrating over
the entire complex $\alpha$ plane, it follows that there exists a
resolution of identity in the form
\begin{equation}
\int d\mu_{f} (\alpha)\,\,|\alpha\rangle_{f}\,_{f}\langle\alpha|
= {\sf I},\qquad d\mu_{f} (\alpha)\,=\,\frac{\rho\,d\rho\,d\theta}
{\pi}\,g_{f}(\rho^{2}),
\label{eq:fresid}
\end{equation}
where the measure $g_{f} (\rho^{2})$ obeys an infinite number of
moment relations:
\begin{equation}
2\,\int_{0}^{\infty}\,d\rho\,\rho^{2 n + 1}\,g_{f} (\rho^{2})\,
(N_{f}(\rho^{2}))^{2}\,=\,n!\,\Gamma (n + 2 k)\,([f(n + k)]!)^{2}
\quad \forall n = 0, 1,2,\cdots.
\label{eq:fmeasure}
\end{equation}
Consequently the measure $g_{f} (\rho^{2})$ may be explicitly
obtained in terms of the inverse Mellin transform as
\begin{equation}
g_{f} (\rho^{2}) = \frac{1}{2 \pi i\,(N_{f} (\rho^{2}))^{2}}
\int_{c - i\infty}^{c + i \infty} ds\,\rho^{- 2 s}\,\Gamma (s)\,
\Gamma (2 k + s - 1)\,([f(k + s - 1)]!)^{2},
\label{eq:Melinv}
\end{equation}
provided $[f(k + n)]!$ may be continued suitably over the region of
integration in the $s$ plane. The construction (\ref{eq:snbgcs})
of the generalized nonlinear Barut-Girardello coherent states
for the $su_{f}(1, 1)$ algebra previously appeared in
\cite{WSP00}. But the resolution of unity in the form of an ordinary
integral over the complex plane, and the explicit evaluation of the
corresponding measure function for the $su_{q}(1, 1)$ algebra
discussed below were not, to our knowledge, obtained earlier.

\par

After setting the general formalism, here we briefly discuss the
single-node Barut-Girardello coherent state of the $q$-deformed
$su_{q}(1,1)$ algebra \cite{KS97}. This state was previously
studied in \cite{SCSJ94}. The commutation rules and the hermiticity
restrictions for the generators of the $su_{q}(1,1)$ algebra read
\begin{equation}
[\cK_{0}, \cK_{\pm}] = \pm \cK_{\pm},\quad
[\cK_{-}, \cK_{+}] = [2\,\cK_{0}]_{q},\quad
\cK_\pm^{\dagger} = \cK_\mp,\quad \cK_{0}^{\dagger} = \cK_{0},
\label{eq:qalg}
\end{equation}
where $[x]_{q} =(q^{x} - q^{-x})/(q - q^{-1})$. For the purpose of
our work we treat the deformation parameter $q$ as a real number
satisfying $0 < q < 1$. The Hopf structure of the algebra
introduces a noncocommutative coproduct map of the generators,
given by
\begin{equation}
\btf (\cK_{0}) = 1 \otimes \cK_{0} + \cK_{0} \otimes 1,\qquad
\btf (\cK_{\pm}) = q^{\cK_{0}} \otimes \cK_{\pm} +
\cK_{\pm} \otimes q^{- \cK_{0}}.
\label{eq:qcoprod}
\end{equation}
A well-known property of the Hopf algebra states that the coproduct
map is a homomorphism of the algebra, given here in (\ref{eq:qalg}).
This will be later used in Sec. 3 to obtain, in the context of a
bipartite composite system, the Barut-Girardello coherent states
for the $q$-deformed $su_{q}(1,1)$ algebra. In the present purpose
of obtaining the Barut-Girardello coherent states of a single-node
system, the coproduct map is, however, not relevant.
\par

Using the well-known Curtright-Zachos \cite{CZ90} map, valid for
a generic $q$, the representations of the $q$-deformed $su_{q}(1,1)$
algebra may be obtained {\it via} the unitary representation
(\ref{eq:clrep}) of the classical $su(1,1)$ algebra. The mapping
function, introduced for the general deformation in
(\ref{eq:fsumap}), here reads
\begin{equation}
f(K_{0}) = \sqrt{\frac{[K_{0} - k]_{q}\,[K_{0} + k - 1]_{q}}
{(K_{0} - k)\,(K_{0} + k - 1)}}.
\label{eq:fsuq}
\end{equation}
The normalized single-node Barut-Girardello coherent state for the
$q$-deformed $su_{q}(1,1)$ algebra is defined as
\begin{equation}
\cK_{-}\,|\alpha, k\rangle_{q} = \alpha\,|\alpha, k\rangle_{q}
\label{eq:qsndef}
\end{equation}
and is explicitly given by
\begin{equation}
|\alpha, k\rangle_{q} = \frac{|\alpha|^{k - \half}}
{\sqrt{I_{2 k -1}^{(q)} (2 |\alpha|)}}\,\sum_{n = 0}^{\infty}
\frac{\alpha^{n}}{\sqrt{[n]_{q}!\,[n+2 k - 1]_{q}!}}\,
|n, k\rangle_{q},\qquad[n]_{q}! = \prod_{j = 1}^{n}[j]_{q}!,
\label{eq:qsnco}
\end{equation}
where the $q$-deformed modified Bessel function reads
\begin{equation}
I_{m}^{(q)} (2 z) = \sum_{n = 0}^{\infty}\frac{z^{m + 2 n}}
{[n]_{q}!\,[m + n]_{q}!}.
\label{eq:qBess}
\end{equation}
The completeness of the $q$-deformed Barut-Girardello coherent
states obtained above may be demonstrated, and the corresponding
measure $g_{q}(\rho^{2})$ defined {\it \`{a} la} (\ref{eq:fresid})
and (\ref{eq:Melinv}) may be explicitly obtained. We will only
present the result here \cite{CV03}:
\begin{eqnarray}
&&g_{q}(\rho^{2}) = \half \, I_{\nu}^{(q)}(2 \rho)
\left[ \frac{q^{2} - 1}{q \ln q}\,\sum_{l = 0}^{\nu -1} (-1)^{l}
\,\frac{[\nu - l -1]_{q}!}{[l]_{q}!}\,\rho^{2 l - \nu}
\right.\nonumber\\
&&\phantom {g_{q}(\rho^{2}) =} + (- 1)^{\nu + 1}\,
\frac{(1 - q^{2})^{2}}{q^{2} (\ln q)^{2}}\,\sum_{l = 0}^{\infty}
\,\frac{1}{[l]_{q}!\,[l + \nu]_{q}!}\,\Big( \ln \rho - \half\,
\psi_{q^{2}}(l + 1) - \half \,\psi_{q^{2}}(l + \nu + 1)\nonumber\\
&&\phantom {g_{q}(\rho^{2}) =}\left. +\half \, (2 l + \nu - 3)\,
\ln q \Big)\,\rho^{2 l + \nu} \right],
\label{eq:qmeasure}
\end{eqnarray}
where $\nu = 2 k - 1$. To obtain the analytical continuation of
$q$-factorial $[n]_{q}!$ we have used the $q$-Gamma function
\cite{GR91} defined as
\begin{equation}
\Gamma_{q}(z) = (1 - q)^{1 - z}\,\frac{(q; q)_{\infty}}
{(q^{z}; q)_{\infty}},\quad \psi_{q}(z) = \frac{d}{dz}\,
\hbox{ln}\,\Gamma_{q}(z),\quad (a; q)_{n} = \sum_{j = 1}^{n}\,
(1 - a\,q^{j - 1}).
\label{eq:qGamma}
\end{equation}
To our knowledge the above measure relating to the completeness
of the single-node Barut-Girardello coherent state for the
$su_{q}(1, 1)$ algebra has not appeared elsewhere. The relevant
classical measure \cite{BG71} in the $q \rightarrow 1$ limit is
readily obtained from (\ref{eq:qmeasure}) as
\begin{equation}
g_{q}(\rho^{2}) \longrightarrow g(\rho^{2}) = 2 I_{\nu}(2 \rho)\,
K_{\nu}(2 \rho),
\label{eq:clmeasure}
\end{equation}
where $K_{\nu}(2 \rho)$ is the modified Bessel function of the
second kind given by
\begin{eqnarray}
&&K_{\nu}(2 \rho) = \half\,\sum_{l= 0}^{\nu -1}\,(- 1)^{l}\,
\frac{(\nu - l - 1)!}{l!}\,\rho^{2 l - \nu}\nonumber\\
&&{\phantom {K_{\nu}(2 \rho) =}} + (- 1)^{\nu + 1}\,
\sum_{l = 0}^{\infty}\frac{1}{l! (l + \nu)!}\,
\Big( {\hbox {ln}} \rho -\half \psi(l + 1) -
\half \psi(l + \nu + 1) \Big)\,\rho^{2 l + \nu},
\label{eq:Kexp}
\end{eqnarray}
where $\psi (z) = (\hbox{ln}\,\Gamma (z))^{\prime}$.

\sect{Bipartite composite system}

Our objective in this section is to construct normalized
Barut-Girardello coherent state for $q$-deformed $su_{q}(1,1)$
algebra in the case of a bipartite composite system. As a
benchmark we first consider this problem for the classical
$su(1,1)$ algebra, where the relevant state is defined by
\begin{equation}
\btf(K_{-})\,|\alpha; k_{1}, k_{2}\rangle =
\alpha\,|\alpha; k_{1}, k_{2}\rangle,\qquad
(2 k_{1}, 2 k_{2}) = \pm 1, \pm 2,\cdots.
\label{eq:clbico}
\end{equation}
The classical coproduct property (\ref{eq:clcoprod}) immediately
provides the factorized form:
\begin{equation}
|\alpha; k_{1}, k_{2}\rangle = |\alpha_{1}, k_{1}\rangle
\otimes|\alpha_{2}, k_{2}\rangle,\qquad
(\alpha_{1}, \alpha_{2}) \in {\mathbb C},
\label{eq:clfact}
\end{equation}
+where
\begin{equation}
\alpha = \alpha_{1} + \alpha_{2}.
\label{eq:alpsum}
\end{equation}
Expanding the state (\ref{eq:clfact}) in a tensored basis
\begin{equation}
|\alpha; k_{1}, k_{2}\rangle = \sum_{n_{1} = 0}^{\infty}\,
\sum_{n_{2} = 0}^{\infty}\,c_{n_{1}, n_{2}}\,|n_{1}, k_{1}\rangle
\otimes |n_{2}, k_{2}\rangle,
\label{eq:tnclexp}
\end{equation}
and using the classical coproduct structure (\ref{eq:clcoprod})
we obtain a recurrence relation
\begin{equation}
\sqrt{(n_{1} + 1)\,(n_{1} + 2 k_{1})}\,c_{n_{1} + 1, n_{2}} +
\sqrt{(n_{2} + 1)\,(n_{2} + 2 k_{2})}\,c_{n_{1}, n_{2} + 1}
= (\alpha_{1} + \alpha_{2})\,c_{n_{1}, n_{2}}.
\label{eq:clrec}
\end{equation}
The solution of the above double-indexed recurrence relation
may be given, of course, as
\begin{equation}
c_{n_{1}, n_{2}} = N_{1}\,N_{2}\,\frac{\alpha_{1}^{n_{1}}\,
\alpha_{2}^{n_{2}}}{\sqrt{n_{1}!\,n_{2}!\,\Gamma(n_{1} + 2 k_{1})
\,\Gamma(n_{2} + 2 k_{2})}},
\label{eq:clrecsol}
\end{equation}
where the normalization factors $N_{i}$ for $i \in (1, 2)$ may be
directly read from (\ref{eq:cslim}) as
\begin{equation}
N_{i} = \frac{|\alpha_{i}|^{k_{i} - \half}}
{\sqrt{I_{2 k_{i} - 1}(2 |\alpha_{i}|)}}.
\label{eq:clnorm}
\end{equation}
We have recapitulated the above facts for the purpose of easy
comparison with our following construction of a bipartite
Barut-Girardello coherent state for $q$-deformed $su_{q}(1,1)$
algebra.

\par

We proceed by defining the said bipartite coherent state for the
$su_{q}(1,1)$ algebra as an eigenstate of the tensored operator
$\btf({\cK_{-})}$:
\begin{equation}
\btf({\cK_{-})}\,|\alpha; k_{1}, k_{2}\rangle_{q} =  \alpha\,
|\alpha; k_{1}, k_{2}\rangle_{q}.
\label{eq:qcodef}
\end{equation}
The state $|\alpha; k_{1}, k_{2}\rangle_{q}$ may again be expanded
{\it {\`a} la} (\ref{eq:tnclexp}) as
\begin{equation}
|\alpha; k_{1}, k_{2}\rangle_{q} = \sum_{n_{1} = 0}^{\infty}\,
\sum_{n_{2} = 0}^{\infty}\,c_{n_{1}, n_{2}}^{(q)}\,
|n_{1}, k_{1}\rangle \otimes |n_{2}, k_{2}\rangle.
\label{eq:qexp}
\end{equation}
The noncocommutative coproduct structure (\ref{eq:qcoprod}) in
conjunction with the mapping function (\ref{eq:fsuq}) now yield
a double-indexed recurrence relation for the above coefficients as
\begin{equation}
q^{n_{1} + k_{1}}\,\sqrt{[n_{1} + 1]_{q}\,[n_{1} + 2 k_{1}]_{q}}\,
c_{n_{1} + 1, n_{2}}^{(q)} + q^{-n_{2} - k_{2}}\,
\sqrt{[n_{2} + 1]_{q}\,[n_{2} + 2 k_{2}]_{q}}\,
c_{n_{1}, n_{2} + 1}^{(q)} = \alpha\,c_{n_{1}, n_{2}}^{(q)}.
\label{eq:qtnrec}
\end{equation}
In the followings we outline a procedure employed here for solving
the above recurrence relation.

\par

We notice that as $q \rightarrow 1$, the deformed recurrence
relation (\ref{eq:qtnrec}) reduces to its classical analogue
(\ref{eq:clrec}), provided the constraint (\ref{eq:alpsum}) is
maintained. In obtaining the solution of the quantized recursion
relation (\ref{eq:qtnrec}) in the presence of the constraint
(\ref{eq:alpsum}), we mimic the classical solution
(\ref{eq:clrecsol}) and consider the following ansatz:
\begin{equation}
c_{n_{1}, n_{2}}^{(q)} = \frac{\alpha_{1}^{n_{1}}\,
\alpha_{2}^{n_{2}}}{\sqrt{[n_{1}]_{q}!\,[n_{2}]_{q}!\,
[n_{1} + 2 k_{1} - 1]_{q}!\,[n_{2} + 2 k_{2} - 1]_{q}!}}\,
g_{n_{1}, n_{2}},
\label{eq:qrecsol}
\end{equation}
where the $q$-dependent coefficients $g_{n_{1}, n_{2}}$ are yet to
be determined. In order to stay close to the classical solution,
we retain the additive property (\ref{eq:alpsum}) for an arbitrary
value of the deformation parameter $q$. For the following
construction the complex parameters $(\alpha_{i}|\,i = (1, 2))$,
while being subjected to the constraint (\ref{eq:alpsum}), are
otherwise arbitrary. Inserting the ansatz (\ref{eq:qrecsol}) in the
recurrence relation (\ref{eq:qtnrec}), we get a simpler recurrence
relation satisfied by the coefficients $g_{n_{1}, n_{2}}$:
\begin{equation}
\alpha_{2}\,q^{n_{1} + k_{1}}\,g_{n_{1}, n_{2} + 1} +
\alpha_{1}\,q^{- n_{2} - k_{2}}\,g_{n_{1} + 1, n_{2}} =
\alpha\,g_{n_{1}, n_{2}}.
\label{eq:grec}
\end{equation}
We impose the limiting condition
\begin{equation}
g_{n_{1}, n_{2}} \rightarrow 1\quad \hbox{as}\quad q \rightarrow 1,
\label{eq:glimit}
\end{equation}
which is consistent with the relation (\ref{eq:alpsum}). Introducing
the parameters
\begin{equation}
\xi = \frac{\alpha_{1}}{\alpha}\,q^{- k_{2}},\quad
\eta = \frac{\alpha_{2}}{\alpha}\,q^{k_{1}}
\label{eq:param}
\end{equation}
and redesignating the indices, we rewrite the recurrence relation
(\ref{eq:grec}) as
\begin{equation}
\eta\,q^{n}\,g_{n, m + 1} + \xi\,q^{- m}\,g_{n + 1, m} = g_{n, m}.
\label{eq:grecfn}
\end{equation}

\par

We now proceed towards solving the above recurrence relation. If we
think of the coefficients $g_{n, m}$ as elements of a matrix, a
little reflection shows that given the elements in the first row,
all other elements can be obtained from (\ref{eq:grecfn})
successively. Accordingly, we assume that elements in the first row
are given as initial conditions:
\begin{equation}
g_{0, m} = d_{m},\qquad m\geq 0.
\label{eq:incon}
\end{equation}
The coefficients $d_{m}$ are {\it arbitrary}, except for the limiting
constraint:
\begin{equation}
d_{m} \rightarrow 1\quad \hbox{as} \quad q \rightarrow 1.
\label{eq:dcon}
\end{equation}
The recurrence relation (\ref{eq:grecfn}) may be systematically used
to completely determine the coefficients $g_{n, m}$ in terms of the
initial distribution $d_{m}$. The emerging pattern suggests the
following ansatz:
\begin{equation}
g_{n, m} = q^{nm}\,\xi^{- n}\,\sum_{k = 0}^{n} (-1)^{k} \eta^{k}
q^{k (k - 1)}\,d_{m + k}\,h_{n, k}(q),
\label{eq:gansatz}
\end{equation}
where the elements $h_{n, k}$ are polynomials in the deformation
parameter $q$, such that
\begin{equation}
h_{0, 0} = 1,\qquad \xi^{-n}\,\sum_{k = 0}^{n} (- 1)^{k}\,
h_{n, k} (q)\Big|_{q \rightarrow 1} = 1.
\label{eq:hlimit}
\end{equation}
Substituting the ansatz (\ref{eq:gansatz}) in the recurrence
relation (\ref{eq:grecfn}) and comparing powers of $\eta$ on both
sides, we get
\begin{equation}
h_{n, 0} = 1,\qquad h_{n, n} = h_{n - 1, n - 1} = \cdots =
h_{0, 0} = 1
\label{eq:hval}
\end{equation}
and the recurrence relation
\begin{equation}
h_{n + 1, k} = h_{n, k} + q^{2 (n - k + 1)}\,h_{n, k - 1}
\quad\hbox{for}\, 1 \leq k \leq n.
\label{eq:hrec}
\end{equation}
A clue to the solution of the recurrence relation (\ref{eq:hrec})
is provided by its classical limit:
\begin{equation}
h_{n + 1, k}^{(q \rightarrow 1)} = h_{n, k}^{(q \rightarrow 1)}
+ h_{n, k - 1}^{(q \rightarrow 1)},
\label{eq:clhrec}
\end{equation}
whose well-known solution reads
\begin{equation}
h_{n, k}^{(q \rightarrow 1)} = \left( \begin{array}{c}n\\k
\end{array} \right).
\label{eq:clhsol}
\end{equation}
This strongly suggests that in the $q$-deformed case $h_{n, k}$
involves $q$-binomial coefficients. In view of the classical
solution (\ref{eq:clhsol}), we try the following ansatz:
\begin{equation}
h_{n, k} = \left\{ \begin{array}{c}n\\k \end{array} \right\}_{q^2}
\equiv \frac{(q^{2};q^{2})_{n}}{(q^{2};q^{2})_{k}\;
(q^{2};q^{2})_{n - k}}.
\label{eq:qhsol}
\end{equation}
This indeed solves the recurrence relation (\ref{eq:hrec}) and
yields the correct classical limit (\ref{eq:clhsol}). The solution
of the recurrence relation (\ref{eq:grecfn}) may now be constructed
{\it via} (\ref{eq:gansatz}) as
\begin{equation}
g_{n, m} = q^{nm}\,\xi^{- n}\,\sum_{k = 0}^{n} (-1)^{k} \eta^{k}
q^{k (k - 1)}\,d_{m + k}\,\left\{ \begin{array}{c}n\\k\end{array}
\right\}_{q^2}.
\label{eq:qgsol}
\end{equation}
The condition (\ref{eq:dcon}) readily yields the limiting value:
${\bf{\big(}}g_{n, m}{\bf{\big)}}_{q \rightarrow 1} = 1$. For special
choices of the boundary coefficients $d_{m}$ the right hand side of
(\ref{eq:qgsol}) may be expressed in a closed form. For instance,
if we choose
\begin{equation}
d_{m} = \delta^{m},\qquad \delta \rightarrow 1\quad\hbox{as}\quad
q \rightarrow 1,
\label{eq:boundgs}
\end{equation}
we obtain
\begin{equation}
g_{n, m} = q^{nm}\,\delta^{m}\,\xi^{- n}\,(\delta\,\eta;\,
q^{2})_{n} = q^{nm}\,\delta^{m}\,\xi^{- n}\,(1 - q^{2n}\,\delta\,
\eta)_{q^{2}}^{n},
\label{eq:gclosed}
\end{equation}
where the $q$-binomial sum \cite{GR91} is expressed as
\begin{equation}
(x;q)_{n} = \sum_{j = 0}^{n}\, (- 1)^{j} q^{j (j + 1)/2}
\left\{ \begin{array}{c}n\\j \end{array} \right\}_{q^{2}}\, x^{j}.
\label{eq:qbin}
\end{equation}
In the second equality in (\ref{eq:gclosed}) we have used the
notation $(1 - z)^{n}_{q} \equiv (q^{- n}\,z;q)_{n}$.
The norm of the bipartite $q$-deformed Barut-Girardello coherent
state introduced in (\ref{eq:qcodef}) is now readily obtained
{\it via} (\ref{eq:qexp}), (\ref{eq:qrecsol}) and
(\ref{eq:qgsol}). Here, for the purpose of simplicity, we
explicitly consider the boundary condition (\ref{eq:boundgs}).
Using the closed form expression (\ref{eq:gclosed}) of the
coefficients $g_{n, m}$, we now obtain the said norm as
\begin{eqnarray}
{\cal N}^{- 2} &{\equiv}&
_{q}\langle \alpha; k_{1}, k_{2}|\alpha; k_{1}, k_{2}\rangle_{q}
= \sum_{n = 0}^{\infty}\,\sum_{m = 0}^{\infty}
|c_{n, m}^{(q)}|^{2}\nonumber\\
\phantom {{\cal N}^{2}} &=& |\delta \alpha_{2}|^{1 - 2 k_{2}}\,
\sum_{n = 0}^{\infty}\frac{q^{n}\,
I_{2 k_{2} - 1}^{(q)}(2\,q^{n}\,|\delta \alpha_{2}|)}
{[n]_{q}!\,[n + 2 k_{1} - 1]_{q}!}\;\big| \alpha^{n}\,
\big(1 - q^{2 n}\,\delta \eta\big)_{q^{2}}^{n} \big|^{2},
\label{eq:qnorm}
\end{eqnarray}
where the modified $q$-Bessel function $I_{m}^{(q)}(2 z)$ is given
in (\ref{eq:qBess}). As the norm has been expressed above as
{\em single-indexed} series sum, its convergence in the domain
$0 < q < 1$ may be tested in a straight-forward way. In the
$q \rightarrow 1$ limit, the norm reduces to its classical value
${\cal N} = N_{1} N_{2}$, where the normalization constants
$N_{i}\:(i = (1, 2))$ are given by (\ref{eq:clnorm}). For the
domain $q > 1$, we may replace the coefficients $g_{n_{1}, n_{2}}$
in (\ref{eq:qrecsol}) by $\hat{g}_{n_{1}, n_{2}} =
g_{n_{2}, n_{1}}(\xi \rightarrow \eta, \eta \rightarrow \xi,
q \rightarrow q^{-1})$ and thereby obtain a finite normed
$q$-deformed bipartite coherent state. This possibility arises on
account of the `crossing symmetry' of the recurrence relation
(\ref{eq:grecfn}) which implies that if $g_{n, m}(\xi, \eta, q)$
is a solution of the said equation, then so is
$g_{m, n}(\eta, \xi, q^{- 1})$.

\par

Combining the above derivation we now present the promised
bipartite normalized Barut-Girardello coherent state for the
quantized $su_{q}(1, 1)$ algebra:
\begin{eqnarray}
&&|\alpha; k_{1}, k_{2}\rangle_{q} = {\cN}\,\sum_{n_{1} = 0}^{\infty}
\,\sum_{n_{2} = 0}^{\infty}\,q^{n_{1} n_{2}} \delta^{n_{2}}
\xi^{- n_{1}} (\delta \eta; q^{2})_{n_{1}}\,\times\nonumber\\
&&\phantom{|\alpha; k_{1}, k_{2}\rangle_{q} =} \times
\Big(\prod_{i = 1}^{2}\,\frac{\alpha_{i}^{n_{i}}}
{\sqrt{[n_{i}]_{q}! [n_{i} + 2k_{i} - 1]_{q}!}}\Big)\,
|n_{1}, k_{1}\rangle \otimes |n_{2}, k_{2}\rangle,
\label{eq:qbBGcost}
\end{eqnarray}
where we have chosen the boundary condition (\ref{eq:boundgs})
for the purpose of simplicity. It is evident that in the classical
$q \rightarrow 1$ limit the above state reduces to the factorized
form (\ref{eq:clfact}).

\sect{Conclusion}

In a bipartite composite system we constructed normalizable
Barut-Girardello coherent state for a quantized $su_{q}(1, 1)$
algebra. Its most remarkable property, as evidenced in
(\ref{eq:qbBGcost}), is the existence of a natural entangled
structure for a nonclassical value of the deformation parameter
$(q \neq 1)$. It is evident from the fact that the summand in
(\ref{eq:qbBGcost}) includes a term $q^{n_{1} n_{2}}$ which
forbids factorization of the relevant coherent state of the
composite system into quantum states of single-node subsystems.
In the classical $q \rightarrow 1$ limit, the entanglement of the
state (\ref{eq:qbBGcost}) disappears as it, in that limit, gets
factorized to the form (\ref{eq:clfact}). Another aspect of the
present derivation is that a {\it one parameter class} of
deformed coherent states with arbitrarily {\it distinct} choices of
boundary values of $d_{m}$, subject to the limiting constraint
(\ref{eq:dcon}), goes to the unique classical limit
(\ref{eq:clfact}) as $q \rightarrow 1$.

\par

The underlying reason of the present structure of entanglement
is the noncocommutativity of the coproduct structure of the quantum
algebras. Therefore entangled structure of the coherent states of
composite systems with quantized symmetries is likely to be a
generic feature. Exploitation of these entangled states obtained
here in the context of quantum teleportation \cite{BBCJPW93} and
entanglement swapping \cite{ZZHE93} is under study. Lastly we
remark that the general nonlinear deformed $su_{f}(1, 1)$ algebra
(\ref{eq:fsualg}) may also be used to obtain entangled states of
bipartite composite systems, as an induced coproduct structure may
be suitably imparted to this algebra. The nature of these
entanglements may be quite distinct from the one presented here.

\bigskip
\noindent {\Large{\bf Acknowledgements}}
\medskip
\noindent
{}

Part of the work was done when one of us (RC) visited
Institute of Mathematical Sciences, Chennai, 600 113, India.
He is partially supported by the grant DAE/2001/37/12/BRNS,
Government of India.

\medskip

\bibliographystyle{amsplain}

\end{document}